\documentclass[12pt,reqno]{amsart}
\usepackage{amsfonts}
\usepackage{amssymb}
\numberwithin{equation}{section}
\theoremstyle{plain}
 \newtheorem{thm}{Theorem}[section]

 \newtheorem{prop}[thm]{Proposition}
\theoremstyle{definition}

 \newtheorem{rem}[thm]{Remark}
\newcommand{\eqd}{\overset{\mathrm d}{=}}
\newcommand{\les}{\leqslant}
\newcommand{\ges}{\geqslant}

\begin{document}
\setlength{\parindent}{1.8pc}
\begin{center}
{\bf SOME PROPERTIES OF EXPONENTIAL INTEGRALS\\
OF L\'EVY PROCESSES AND EXAMPLES}\\[5mm]
HITOSHI KONDO\\
{\it Department of Mathematics, Keio University, Hiyoshi, Yokohama, 223-8522 
Japan}\\
email: jin\_kondo@2004.jukuin.keio.ac.jp\\[3mm]
MAKOTO MAEJIMA\\
{\it Department of Mathematics, Keio University, Hiyoshi, Yokohama, 223-8522 
Japan}\\
email: maejima@math.keio.ac.jp\\[3mm]
KEN-ITI SATO\\
{\it Hachiman-yama 1101-5-103, Tenpaku-ku, Nagoya, 468-0074 Japan}\\
email: ken-iti.sato@nifty.ne.jp
\end{center}


\vskip 5mm
{\small {\it Abstract}.\quad
The improper stochastic integral $Z=\int_0^{\infty-}\exp(-X_{s-})dY_s$ is 
studied, where $\{ (X_t ,Y_t) , t \geqslant 0 \}$ is a L\'evy process on 
$\mathbb R ^{1+d}$ with $\{X_t \}$ and $\{Y_t \}$ being $\mathbb R$-valued 
and  $\mathbb R ^d$-valued, respectively. The condition for existence and 
finiteness of $Z$ is given and then the law $\mathcal L(Z)$ of $Z$
is considered. Some sufficient conditions for $\mathcal L(Z)$ to be 
selfdecomposable and some sufficient conditions for $\mathcal L(Z)$ to be 
non-selfdecomposable but semi-selfdecomposable are given. Attention is 
paid to the case where $d=1$, $\{X_t\}$ is a Poisson process, and 
$\{X_t\}$ and $\{Y_t\}$ are independent.  An example of $Z$ of type $G$ 
with selfdecomposable mixing distribution is given.}

\vskip 10mm
%%%%%%%%%%%%%%%%%%%%% Section 1 %%%%%%%%%%%%%%%%%%%%%%%%%%%%%%%%

\section{Introduction}

Let $\{ (\xi _t ,\eta _t) , t \ges 0 \}$ 
be a L\'evy process on $\mathbb R ^{2}$.
The generalized Ornstein-Uhlenbeck process $\{V_t , t \ges 0 \}$ on $\mathbb R$ 
based on $\{ (\xi _t ,\eta _t) , t \ges 0 \}$ with initial condition $V_0$
is defined as
\begin{equation}\label{1.1}
V_t = e^{-\xi _t }\left( V_0 + \int _{0}^{t} e^{\xi _{s-}} {d}\eta _s \right) , 
\quad t\ges 0,
\end{equation}
where $V_0$ is a random variable independent of 
$\{(\xi _t ,\eta _t) , t\ges 0\}$.
This process has recently been well-studied by 
Carmona, Petit, and Yor \cite{CPY97}, \cite{CPY01}, Erickson and Maller \cite{EM05}, and 
Lindner and Maller \cite{LM05}.

Lindner and Maller \cite{LM05} find that the generalized Ornstein-Uhlenbeck 
process $\{V_t , t \ges 0 \}$ based on $\{ (\xi _t ,\eta _t) , t \ges 0 \}$ 
 turns out to be
 a stationary process with a suitable choice of $V_0$ if and only if 
\begin{equation}\label{1.2}
P\left(
\int _{0}^{\infty -} e^{-\xi_{s-}}dL_s \text{ exists and is finite}\right)=1,
\end{equation}
where 
\begin{equation}\label{1.2a}
\int _{0}^{\infty -} e^{-\xi_{s-}}dL_s =
\lim _{t \rightarrow \infty} \int _{0}^{t} e^{-\xi_{s-}}dL_s 
\end{equation}
and $\{ (\xi _t ,L_t) , t \ges 0 \}$ 
is a L\'evy process on $\mathbb R^2$ defined by
\begin{equation}\label{1.3}
L_t=\eta_t+\sum_{0<s\les t}(e^{-(\xi_s-\xi_{s-})}-1)(\eta_s-\eta_{s-})-ta_{\xi,
\eta}^{1,2}
\end{equation}
with $\left(a_{\xi,\eta}^{j,k}\right)_{j,k=1,2}$ being the Gaussian 
covariance matrix in the L\'evy--Khintchine triplet of the process
$\{(\xi _t ,\eta _t)\}$.
Moreover, if the condition \eqref{1.2} is satisfied, then
the choice of $V_0$ which makes $\{V_t\}$ stationary is unique in law and
\begin{equation}\label{1.4}
\mathcal L(V_0)=\mathcal L\left(\int _{0}^{\infty -} e^{-\xi_{s-}}dL_s \right).
\end{equation}
Here $\mathcal L$ stands for \lq\lq the distribution of".
If $\{ \xi _t, t \ges 0 \}$ and $\{ \eta _t , t \ges 0 \}$ are independent,
then $P(L_t=\eta_t\text{ for all }t)=1$.

Keeping in mind the results in the preceding paragraph,
we study in this paper the exponential integral 
$\int _{0}^{\infty -} e^{-X_{s-}}dY_s$, where $\{ (X_t ,Y_t) , t \ges 0 \}$ 
is a L\'evy process on $\mathbb R ^{1+d}$ with $\{X_t \}$ and $\{Y_t \}$
being $\mathbb R$-valued and  $\mathbb R ^d$-valued, respectively. 
In Section 2 the existence conditions for this integral are given.
They complement a theorem for $d=1$ of Erickson and Maller \cite{EM05}.
Then, in Section 3, some
properties of 
\begin{equation}\label{1.5}
\mu=\mathcal L\left(\int _{0}^{\infty -} e^{-X_{s-}}dY_s \right)
\end{equation}
are studied.
A sufficient condition for $\mu$ to be a selfdecomposable distribution on $\mathbb R^d$ 
is given as in Bertoin, Lindner, and Maller \cite{BLM05}.  
Further we give a sufficient condition for $\mu$ 
not to be selfdecomposable.  Recall that, in the case where $X_t=t$, $t\ges0$, and 
$\{Y_t\}$ is a L\'evy process on $\mathbb R^d$, 
$\mathcal L\left(\int _{0}^{\infty -} e^{-s}dY_s
\right)$ is always selfdecomposable if the integral exists and is finite (see 
e.\,g.\ \cite{S}, Section 17).
In particular, we are interested in the case where 
$\{X_t \}$ and $\{Y_t \}$ are independent and $\{X_t \}$ is a Poisson process;
we will give a sufficient condition for $\mu$ to be 
semi-selfdecomposable and not selfdecomposable and also a sufficient condition 
for $\mu$ to be selfdecomposable. 
In Section 4, we are concerned with $\mu$ of \eqref{1.5} when $\{X_t\}$ is a 
Brownian motion with positive drift on $\mathbb R$, $\{Y_t\}$ is a symmetric 
$\alpha$-stable L\'evy process on $\mathbb R$ with $0<\alpha\les2$, and $\{X_t\}$ and
$\{Y_t\}$ are independent. We will show that in this case
$\mu$ gives a type $G$ distribution
with selfdecomposable mixing distribution, which is related to results in 
Maejima and Niiyama \cite{MNi05} and Aoyama, Maejima, and Rosi\'nski \cite{AMR06}.

%%%%%%%%%%%%%%%%%%%%% Section 2%%%%%%%%%%%%%%%%%%%%%%%%%%%%%%%%

\vskip 5mm
\section{Existence of exponential integrals of L\'evy processes}

Let $\{ (X_t ,Y_t) , t \ges 0 \}$ be a L\'evy process on $\mathbb R ^{1+d}$, where
$\{X_t \}$ is $\mathbb R$-valued and $\{Y_t \}$ is $\mathbb R ^d$-valued. 
We keep this set-up throughout this section.
Let $(a_X,\nu_X,
\gamma_X)$ be the L\'evy-Khintchine triplet of the process $\{X_t\}$ in the
sense that
\[
Ee^{izX_t}=\exp\left[t\left(-\frac12 a_X z^2+i\gamma_X z
+\int_{\mathbb R\setminus\{0\}} (e^{izx}-1-izx 1_{\{|x|\les1\}}(x))\nu_X(dx)
\right)\right]
\]
for $z\in\mathbb R$, where $a_X\ges0$ and $\nu_X$ is the L\'evy measure of $\{X_t\}$. 
Denote
\begin{equation}\label{2.0}
h_X(x)=\gamma_X +\nu_X(\,(1,\infty)\,)+\int_1^x\nu_X(\,(y,\infty)\,)
dy.
\end{equation}
Let $\nu_Y$ be the L\'evy measure of $\{Y_t\}$.  
The following result is a
$d$-dimensional extension of Theorem 2 of Erickson and Maller \cite{EM05}.

\begin{thm}\label{t2.1}
Suppose that there is $c>0$ such that $h_X(x)>0$ for all $x\ges c$ and that
$\{Y_t\}$ is not the zero process. Then
\begin{equation}\label{2.1}
P\left(
\int _{0}^{\infty -} e^{-X_{s-}}dY_s \text{ exists and is finite}\right)=1
\end{equation}
if and only if
\begin{equation}\label{2.2}
\lim_{t\to\infty}X_t=+\infty\text{ a.\,s.\ and }\int_{|y|\ges e^c}
\frac{\log|y|}{h_X(\log|y|)} \nu_Y (dy)<\infty,
\end{equation}
where $|y|$ is the Euclidean norm of $y \in \mathbb R^d$.
\end{thm}

\begin{proof}
First, for $d=1$, this theorem is established in \cite{EM05}.
Second, for $j=1,\ldots,d$, the $j$th coordinate process $\{Y_t^{(j)}, t\ges0\}$
is a L\'evy process on $\mathbb R$ with L\'evy measure 
$\nu_{Y^{(j)}}(B)=\int_{\mathbb R^d} 1_B(y_j)\nu_Y (dy)$ for any Borel set $B$ in 
$\mathbb R$
satisfying $0\not\in B$, where $y=(y_1,\ldots,y_d)$. Third,
the property \eqref{2.1} is equivalent to
\begin{equation}\label{2.1a}
P\left(
\int _{0}^{\infty -} e^{-X_{s-}}dY_s^{(j)} \text{ exists and is finite}\right)=1
\quad\text{for }j=1,\ldots,d.
\end{equation}

Next, we claim that the following \eqref{2.3} and \eqref{2.4}
are equivalent:
\begin{gather}
\int_{|y|>M}\frac{\log |y|}{h_X(\log |y|)} \nu_Y(dy)<\infty\quad
\text{ for some }M\ges e^c,\label{2.3}\\
\int_{\{y\colon |y_j|>M\}}
\frac{\log |y_j|}{h_X(\log |y_j|)} \nu_Y(dy)<\infty,\quad
j=1,\ldots,d,\quad\text{for some }M\ges e^c.\label{2.4}
\end{gather}
Put $f(u)=\log u/h_X(\log u)$ for $u\ges e^c$. This $f(u)$ is not necessarily
increasing for all $u\ges e^c$. We use the words {\em increasing} and 
{\em decreasing} in the wide sense allowing flatness.  But
$f(u)$ is increasing for sufficiently large $u$ ( $>M_0$, say), because, for $x>c$,
\[
\frac{h_X(x)}{x}=\frac{h_X(c)}{x}+\frac{1}{x} \int_c^x n(y)dy
\]
with $n(y)=\nu_X(\,(y,\infty)\,)$ and, with 
$d/dx$ meaning the right derivative, we have
\begin{align*}
&\frac{d}{dx}\left( \frac{1}{x} \int_c^x n(y)
dy\right)=\frac{1}{x^2}\left(-\int_c^x n(y)
dy+xn(x)\right)\\
&\qquad=\frac{1}{x^2}\left(\int_c^x (n(x)-n(y))
dy+cn(x)\right)<0
\end{align*}
for sufficiently large $x$
if $n(c)>0$ (note that $\int_c^x (n(x)-n(y))dy$
is nonpositive and decreasing).
Thus we see that \eqref{2.3} implies \eqref{2.4}.
Indeed, letting $M_1 = M\lor M_0$, we have
\[
\int_{\{y\colon |y_j|>M_1\}} f(|y_j|)\nu_Y(dy)\les \int_{\{y\colon 
|y_j|>M_1\}} f(|y|)\nu_Y(dy)\les 
\int_{|y|>M_1} f(|y|)\nu_Y(dy)<\infty.
\]
In order to show that \eqref{2.4} implies \eqref{2.3}, let $g(x)=h_X(x)$ for 
$x\ges c$ and $=h_X(c)$ for $-\infty< x<c$.
Then $g(x)$ is positive and increasing on $\mathbb R$.
Assume \eqref{2.4}. Let $M_1 = M\lor M_0$. 
Then, using the concavity of $\log(u+1)$ for $u\ges0$,
we have
{\allowdisplaybreaks
\begin{align*}
&\int_{|y|>M_1} f(|y|)\nu_Y(dy)\les\int_{|y|>M_1} f(|y_1|+\cdots+|y_d|)
\nu_Y(dy)\\
&\qquad\les
\int_{|y|>M_1}\frac{\log(|y_1|+\cdots+|y_d|+1)}{h_X(\log(|y_1|+\cdots+|y_d|))}
\nu_Y(dy)\\
&\qquad\les
\sum_{j=1}^d \int_{|y|>M_1}\frac{\log(|y_j|+1)}{h_X(\log(|y_1|+\cdots+|y_d|))}
\nu_Y(dy),\\
&\qquad=\sum_{j=1}^d \int_{|y|>M_1}\frac{\log(|y_j|+1)}{g(\log(|y_1|+\cdots+|y_d|))}
\nu_Y(dy)\\
&\qquad\les\sum_{j=1}^d \int_{|y|>M_1}\frac{\log(|y_j|+1)}{g(\log(|y_j|))}
\nu_{\eta}(dy)\\
&\qquad\les\sum_{j=1}^d \left(\int_{|y_j|>M_1}
\frac{\log(|y_j|+1)}{g(\log(|y_j|))}
\nu_Y(dy)+\int_{|y_j|\les M_1,\,|y|>M_1}\frac{\log(|y_j|+1)}{g(\log(|y_j|))}
\nu_Y(dy)\right).
\end{align*}
The first integral in each summand is finite due to} \eqref{2.4} 
 and the second integral is also finite because the integrand is bounded.
This finishes the proof of equivalence of \eqref{2.3} and \eqref{2.4}.

Now assume that \eqref{2.2} holds. Then \eqref{2.4} holds. Hence, by the theorem
for $d=1$, $\int _{0}^{\infty -} e^{-X_{s-}}dY_s^{(j)}$ exists and is finite 
a.\,s.\ for all $j$ such that $\{Y_t^{(j)}\}$ is not the zero process.
For $j$ such that $\{Y_t^{(j)}\}$ is the zero process, we have
$\int _{0}^{\infty -} e^{-X_{s-}}dY_s^{(j)}=0$. Hence \eqref{2.1a} holds, that is,
\eqref{2.1} holds.

Conversely, assume that \eqref{2.1} holds.
Let
\[
I_j=\int_{\{y\colon |y_j|\ges e^c\}}
\frac{\log |y_j|}{h_X(\log |y_j|)} \nu_Y(dy).
\]
Since $\{Y_t\}$ is not the zero process, $\{Y_t^{(j)}\}$ is not the zero process
for some $j$. Hence, by the theorem for $d=1$, $\lim_{t\to\infty}X_t=+\infty$ 
a.\,s.\ and $I_j<\infty$ for such $j$. For $j$ such that $\{Y_t^{(j)}\}$ is the 
zero process, $\nu_{Y^{(j)}}=0$ and $I_j=0$.  Hence we have \eqref{2.4} and thus
\eqref{2.2} holds due to the equivalence of \eqref{2.3} and \eqref{2.4}.
\end{proof}

\begin{rem}\label{r2.1}
(i) Suppose that $\{X_t\}$ satisfies $0<EX_1<\infty$. Then $\lim_{t\to\infty}X_t=
+\infty$ a.\,s.\ and $h_X(x)$ is positive and bounded for large $x$.  Thus
\eqref{2.1} holds if and only if
\begin{equation}\label{2.5}
\int_{\mathbb R^d} \log^+ |y|\,\nu_Y(dy)<\infty.
\end{equation}
Here $\log^+ u=0\lor\log u$.
For $d=1$ this is mentioned in the comments following Theorem 2 of \cite{EM05}.

(ii) As is pointed out in Theorem 5.8 of Sato \cite{S01a}, $\lim_{t\to\infty}X_t
=+\infty$ a.\,s.\ if and only if one of the following (a) and (b) holds:\\
(a)  $E(X_1\land 0)>-\infty$ and $0<EX_1\les +\infty$;\\
(b)  $E(X_1\land 0)=-\infty$, $E(X_1\lor 0)=+\infty$, and
\begin{equation}\label{2.6}
\int_{(-\infty,-2)}|x|\left(\int_1^{|x|} \nu_X(\,(y,\infty)\,)dy\right)^{-1}
\nu_X(dx)<\infty.
\end{equation}
In other words, $\lim_{t\to\infty}X_t=+\infty$ a.\,s.\ if and only if one of 
the following (a$'$) and (b$'$) holds:\\
(a$'$)  $E|X_1|<\infty$ and $EX_1>0$;\\
(b$'$)  $\int_1^{\infty} \nu_X(\,(y,\infty)\,)dy=\infty$ and \eqref{2.6} holds.\\
See also Doney and Maller \cite{DM02}.

(iii) If $\lim_{t\to\infty}X_t=+\infty$ a.\,s., then $h_X(x)>0$ for all large
$x$, as is explained in \cite{EM05} after their Theorem 2.
\end{rem}

When $\{X_t\}$ and $\{Y_t\}$ are
independent, the result in Remark \ref{r2.1} (i) can be extended to more general
exponential integrals of L\'evy processes.

\begin{thm}\label{t2.2}
Suppose that $\{X_t\}$ and $\{Y_t\}$ are
independent and that $0<EX_1<\infty$. Let $\alpha>0$. Then
\begin{equation}\label{2.7}
P\left(
\int _{0}^{\infty -} e^{-(X_{s-})^{\alpha}}dY_s \text{ exists and is finite}\right)=1
\end{equation}
if and only if
\begin{equation}\label{2.8}
\int_{\mathbb R^d} (\log^+ |y|)^{1/\alpha}\nu_Y(dy)<\infty.
\end{equation}
\end{thm}

We use the following result, which is a part of Proposition 4.3 of 
\cite{S05b}.

\begin{prop}\label{p2.1}
Let $f$ be a locally square-integrable nonrandom function on $[0,\infty)$
such that there are positive constants $\alpha$, $c_1$, and $c_2$ satisfying
\begin{equation*}
e^{-c_2 s^{\alpha}}\les f(s)\les e^{-c_1 s^{\alpha}}\quad\text{for all large $s$.}
\end{equation*}
Then 
\[
P\left(
\int _{0}^{\infty -} f(s)dY_s \text{ exists and is finite}\right)=1
\]
if and only if \eqref{2.8} holds.
\end{prop}

\begin{proof}[Proof of Theorem \ref{t2.2}] 
Let $E[X_1]=b$. By assumption, $0<b<\infty$.  By the law of large numbers for 
L\'evy processes (Theorem 36.5 of \cite{S}), we have $\lim_{t\to\infty}X_t/t=b$ 
a.\,s. Hence
\[
P(b/2<\xi_t/t<2b\text{ for all large $t$})=1.
\]
Conditioned by the process $\{X_t\}$, the integral $\int_0^t e^{-(X_{s-})^{\alpha}}dY_s$ 
can be considered as that with $X_s$, $s\ges0$,
 frozen while $Y_s$, $s\ges0$, maintains the same randomness.
This is because the two processes are independent.
Hence we can apply Proposition \ref{p2.1}. Thus, if
\eqref{2.8} holds, then
\begin{align*}
&P\left(
\int _{0}^{\infty -} e^{-(X_{s-})^{\alpha}}dY_s \text{ exists and is finite}\right)\\
&\qquad=
E\left[ P\left(\left.\int _{0}^{\infty -} e^{-(X_{s-})^{\alpha}}dY_s 
\text{ exists and is finite}\;\right|\;\{X_t\}\right)\right]=1.
\end{align*}
Conversely, if \eqref{2.8} does not hold, then
\[
P\left(
\int _{0}^{\infty -} e^{-(X_{s-})^{\alpha}}dY_s \text{ exists and is finite}\right)=0.
\]
Indeed, in the situation of Proposition \ref{p2.1}, we have, 
by Kolmogorov's zero-one law,
\[
P\left(
\int _{0}^{\infty -} f(s)dY_s \text{ exists and is finite}\right)=0
\]
if and only if \eqref{2.8} does not hold. 
\end{proof}

\vskip 5mm
%%%%%%%%%%%%%%%%%%%%% Section 3 %%%%%%%%%%%%%%%%%%%%%%%%%

\section{Properties of the laws of
exponential integrals of L\'evy processes.}

Let $\mu$ be a distribution on $\mathbb R^d$. Denote by $\widehat\mu(z)$, 
$z\in\mathbb R^d$, the
characteristic function of $\mu$. We call $\mu$ selfdecomposable if, for every
$b\in(0,1)$, there is a distribution $\rho_b$ on $\mathbb R^d$ such that
\begin{equation}\label{3.1}
\widehat\mu(z)=\widehat\mu(bz)\widehat\rho_b(z).
\end{equation}
If $\mu$ is selfdecomposable, then $\mu$ is infinitely divisible and $\rho_b$
is uniquely determined and infinitely divisible.
If, for a fixed $b\in(0,1)$, there is an infinitely divisible
 distribution $\rho_b$ on $\mathbb R^d$ satisfying
\eqref{3.1}, then $\mu$ is called $b$-semi-selfdecomposable, or of class 
$L_0(b,\mathbb R^d)$.
If $\mu$ is $b$-semi-selfdecomposable, then $\mu$ is infinitely divisible and $\rho_b$
is uniquely determined.  
If $\mu$ is $b$-semi-selfdecomposable and $\rho_b$
is of class $L_0(b,\mathbb R^d)$, then $\mu$ is called of class $L_1(b,\mathbb R^d)$.
These \lq\lq semi"-concepts were introduced by Maejima and Naito \cite{MNa98}.

We start with a sufficient condition for selfdecomposability of the laws of
exponential integrals of L\'evy processes.

\begin{thm}\label{t3.1}
Suppose that $\{ (X_t ,Y_t) , t \ges 0 \}$
is a L\'evy process on $\mathbb R ^{1+d}$, where
$\{X_t \}$ is $\mathbb R$-valued and $\{Y_t \}$ is $\mathbb R ^d$-valued.
Suppose in addition that $\{X_t \}$ does not have positive jumps and $0<EX_1<+\infty$ 
and that
\begin{equation}\label{3.2}
\int_{\mathbb R^d} \log^+ |y| \nu_Y (dy)<\infty
\end{equation}
for the L\'evy measure $\nu_Y$ of $\{Y_t\}$.  Let
\begin{equation}\label{3.3}
\mu=\mathcal L\left( \int_0^{\infty-} e^{-X_{s-}} dY_s\right).
\end{equation}
Then $\mu$ is selfdecomposable.
\end{thm}

When $d=1$ and $Y_t =t$, the assertion is found in \cite{KLM06}.
When $d=1$, the assertion of this theorem is found 
in the paper \cite{BLM05} with a key idea of the proof. 
This fact was informed personally by 
Alex Lindner to the second author of the present paper
when he was visiting Munich in November, 2005, while the paper \cite{BLM05}
was in preparation.
For $d\ges2$ we do not need a new idea but, for completeness, we give 
a proof of it here.

\begin{proof}[Proof of Theorem \ref{t3.1}]
If $\{Y_t \}$ is the zero process, then the theorem is trivial. Hence we 
assume that $\{Y_t \}$ is not the zero process. 
Under the assumption that $\{X_t \}$ does not have positive jumps, we have that
$\lim_{t\to\infty}X_t=+\infty$ a.\,s.\ if and only if $0<EX_1<+\infty$.
Thus the integral $Z=
\int_0^{\infty-} e^{-X_{s-}} dY_s$ exists and is finite a.\,s.\ by virtue of
Theorem \ref{t2.1} and Remark \ref{r2.1} (i).
Let $c>0$, and define 
$$
T_c = \inf \{ t\colon X_t =c\}.
$$
Since we are assuming that $X_{t}$ does not have positive jumps and that
$0<EX_1<+\infty$, we have $T_c < \infty$ and $X(T_c)=c$ a.\,s.  
 Then we have
\[
Z= \int _{0}^{\infty -} e^{-X_{s-}} {d}Y_s 
= \int _{0}^{T_c}e^{-X_{s-}} {d}Y_s + \int _{T_c}^{\infty -} e^{-X_{s-}} 
{d}Y_s.
\]
Denote by $U_c$ and $V_c$ the first and second integral of the last member.
 We have
\[
V_c = \int _{T_c}^{\infty -} e^{-X({s-})+X({T_c})-X(T_c)} {d}Y_s=e^{-c}Z_c,
\]\\
where 
\[
Z_c=\int _{T_c}^{\infty -} e^{-X({s-})+X(T_c)} {d}Y_s
= \int_0^{\infty -} e^{-(X({T_c+s-})-X({T_c}))}d(Y({T_c+s}) -Y({T_c})).
\]
Since $T_c$ is a stopping time for the process $\{(X_s,Y_s), s\ges0\}$, 
we see that $\{(X(T_c+s)-X(T_c),Y(T_c+s)-Y(T_c)), s\ges0\}$ and
$\{(X_s,Y_s), 0\les s\les T_c\}$ are independent and the former process is
identical in law with $\{(X_s,Y_s), s\ges0\}$ (see Theorem 40.10 of \cite{S}). 
Thus $Z_c$ and $U_c$ are independent and $\mathcal L(Z_c)=\mathcal L(Z)$. 
Since $c$ is arbitrary, it follows 
that the law of $Z$ is selfdecomposable.
\end{proof}

We turn our attention to the case where $\{X_t\}$ is a Poisson process
and $\{X_t\}$ and $\{Y_t\}$ are independent. 
The suggestion of studying this case was personally given by Jan Rosi\'nski
to the authors. In this
case we will show that the law $\mu$ of the exponential integral can be 
selfdecomposable or non-selfdecomposable, depending on the choice of $\{Y_t\}$. 
A measure $\nu$ on $\mathbb R^d$ is called discrete if it is concentrated on some
 countable set $C$, that is, 
 $\nu(\mathbb R^d\setminus C)=0$. 

\begin{thm}\label{t3.2}
Suppose that $\{N_t, t\ges0\}$ is a Poisson process, $\{Y_t\}$ is a L\'evy process
on $\mathbb R^d$,
and $\{N_t\}$ and $\{Y_t\}$ are independent. Suppose that \eqref{3.2} holds. Let
\begin{equation}\label{3.3a}
\mu=\mathcal L\left( \int_0^{\infty-} e^{-N_{s-}} dY_s\right).
\end{equation}
Then the following statements are true.

{\rm(i)} The law $\mu$ is infinitely divisible and, furthermore, 
$e^{-1}$-semi-selfdecomposable.

{\rm(ii)} Suppose that either  $\{Y_t\}$ is a strictly $\alpha$-stable L\'evy process on 
$\mathbb R^d$, $d\ges1$, 
 with $0<\alpha\les2$ or  $\{Y_t\}$ is a Brownian motion with drift with $d=1$.
Then, $\mu$ is selfdecomposable and of class $L_1(e^{-1},\mathbb R^d)$.

{\rm(iii)} Suppose that $d=1$ and $\{Y_t\}$ is integer-valued,
not identically zero. Let
\[
D=\begin{cases}(0,\infty) \quad & \text{if\/ $\{Y_t\}$ is increasing},\\
(-\infty,0) \quad & \text{if\/ $\{Y_t\}$ is decreasing},\\
\mathbb R \quad & \text{if\/ $\{Y_t\}$ is neither increasing nor decreasing}.
\end{cases}
\]
Then $\mu$ is not selfdecomposable and, furthermore, the L\'evy measure 
$\nu_{\mu}$ of $\mu$ is discrete and the set of points with positive 
$\nu_{\mu}$-measure is dense in $D$.
\end{thm}

It is noteworthy that a seemingly pathological L\'evy measure 
appears in a natural way in the assertion (iii).
In relation to the infinite divisibility in (i), we recall that
$\int_0^{\infty-} \exp(-N_{s-}-cs) ds$ does not have an infinitely divisible law
if $c>0$.  This is Samorodnitsky's remark mentioned in \cite{KLM06}.
The integral $\int_0^{\infty-} \exp(-N_{s-}) ds$ is a special case of (ii) with
$\alpha=1$.

\begin{proof}[Proof of Theorem \ref{t3.2}]
(i) Let $Z=\int_0^{\infty-} e^{-N_{s-}} dY_s$.
If $\{Y_t\}$ is the zero process, then $Z=0$.  
If $\{Y_t\}$ is not the zero process, then existence and finiteness of $Z$
follows from Theorem \ref{t2.1}.  Let $T_n=\inf\{s\ges0\colon N_s=n\}$. 
Clearly $T_n$ is finite and tends to infinity as $n\to\infty$ a.\,s. We have
\begin{equation*}
Z=\sum_{n=0}^{\infty}\int_{T_n}^{T_{n+1}} e^{-N_{s-}} dY_s=\sum_{n=0}^{\infty}
e^{-n}(Y(T_{n+1})-Y(T_n)).
\end{equation*}
For each $n$, $T_n$ is a stopping time for $\{(N_s,Y_s)\colon s\ges0\}$. Hence
$\{(N(T_n +s)-N(T_n), Y(T_n +s)-Y(T_n)), s\ges0\}$ and $\{(N_s,Y_s), 0\les s\les
T_n\}$ are independent and the former process is identical in law with 
$\{(N_s,Y_s), s\ges0\}$. It follows that the family
$\{Y(T_{n+1})-Y(T_n), n=0,1,2,\ldots\}$ is independent and identically distributed.
Thus, denoting $W_n= Y(T_{n+1})-Y(T_n)$, we have representation
\begin{equation}\label{3.4}
Z=\sum_{n=0}^{\infty} e^{-n} W_n,
\end{equation}
where $W_0,W_1.\ldots$ are independent and identically distributed
 and $W_n\eqd Y(T_1)$ ( $\eqd$ stands for
\lq\lq has the same law as"). Consequently we have
\begin{equation}\label{3.5}
Z=W_0+e^{-1}Z',
\end{equation}
where $W_0$ and $Z'$ are independent and $Z'\eqd Z$.
The distribution of $W_0$ is infinitely divisible, since $W_0=Y(T_1)\eqd U_1$,
where $\{U_s\}$ is a L\'evy process given by subordination of $\{Y_s\}$ by a 
gamma process.
Here we use our assumption of independence of $\{N_t\}$ and $\{Y_t\}$. Thus
$\mu$ is
$e^{-1}$-semi-selfdecomposable and hence infinitely divisible.
An alternative proof of the infinite divisibility of $\mu$ is to look at 
the representation \eqref{3.4} and to use that $\mathcal L(Y(T_1))$
is infinitely divisible.  

(ii) 
Use the representation \eqref{3.4} with $W_n\eqd U_1$,
where we obtain a L\'evy process $\{U_s\}$ by subordination of $\{Y_s\}$ 
by a gamma process.
Since gamma distributions are selfdecomposable, the results of Sato \cite{S01b}
on inheritance of selfdecomposability in subordination
guarantee that $\mathcal L(U_1)$ is selfdecomposable under our assumption on $\{Y_s\}$.
Hence $\mu$ is selfdecomposable, as selfdecomposability is preserved under 
convolution and convergence.
Further, since selfdecomposability implies $b$-semi-selfdecomposability for each
$b$, \eqref{3.5} shows that $\mu$ is of class $L_1(e^{-1},\mathbb R^d)$.

(iii) 
The process $\{Y_t\}$ is a compound Poisson process on $\mathbb R$ with $\nu_Y$ 
concentrated on the integers (see Corollary 24.6 of \cite{S}). Let us consider
the L\'evy measure $\nu^{(0)}$ of $Y(T_1)$. 
Let $a>0$ be the parameter of the Poisson process $\{N_t\}$.
As in the proofs of (i)
and (ii), $Y(T_1)\eqd U_1$, where 
$\{U_s\}$ is given by subordination of $\{Y_s\}$, by a gamma process
which has L\'evy measure $x^{-1}e^{-ax}dx$. 
Hence, using Theorem 30.1 of \cite{S}, we see that
\begin{equation*}
\nu^{(0)}(B)=\int_0^{\infty}P(Y_s\in B)s^{-1}e^{-as}ds
\end{equation*}
for any Borel set $B$ in $\mathbb R$. Thus $\nu^{(0)}(\mathbb R\setminus\mathbb Z)=0$.

Suppose that $\{Y_t\}$ is not a decreasing process. Then some positive integer
has positive $\nu^{(0)}$-measure. Denote by $p$ the minimum of such positive integers.
Since $\{Y_t\}$ is compound Poisson, $P(Y_s=kp)>0$ for any $s>0$ for 
$k=1,2,\ldots$. Hence $\nu^{(0)}(\{kp\})>0$ for $k=1,2,\ldots$. Therefore,
for each nonnegative
integer $n$, the L\'evy measure $\nu^{(n)}$ of $e^{-n}Y(T_1)$ satisfies
$\nu^{(n)}(\{e^{-n}kp\})>0$ for $k=1,2,\ldots$. 
Clearly, $\nu^{(n)}$ is also discrete.
The representation \eqref{3.4} shows that
\[
\nu_{\mu}=\sum_{n=0}^{\infty} \nu^{(n)}.
\]
Hence, $\nu_{\mu}$ is discrete and
\[
\nu_{\mu}(\{e^{-n}kp\})>0\quad\text{for all $n=0,1,2,\ldots$ and $k=1,2,\ldots$\;\,.} 
\]
Thus the points in $(0,\infty)$ of positive $\nu_{\mu}$-measure are dense in
$(0,\infty)$.

Similarly, if $\{Y_t\}$ is not an increasing process, then the points in 
$(-\infty,0)$ of positive $\nu_{\mu}$-measure are dense in $(-\infty,0)$.
\end{proof}

The following remarks give information on continuity properties of the law $\mu$.
A distribution on $\mathbb R^d$ is called
nondegenerate if its support is not contained in any 
affine subspace of dimension $d-1$. 

\begin{rem}\label{r3.1}
(i) Any nondegenerate selfdecomposable distribution on $\mathbb R^d$ for $d\ges1$
is absolutely continuous (with respect to Lebesgue measure on $\mathbb R^d$)
 although, for $d\ges2$, its L\'evy measure is not necessarily
absolutely continuous. This is proved by Sato \cite{S82} (see also Theorem 27.13
of \cite{S}).  

(ii) Nondegenerate semi-selfdecomposable distributions on $\mathbb R^d$ for $d\ges1$
are absolutely continuous or continuous singular, as Wolfe \cite{W83} proves
(see also Theorem 27.15 of \cite{S}).
\end{rem}

%%%%%%%%%%%%%%%%%%%% section 4 %%%%%%%%%%%%%%%%%%%%%%%%%%%%%%
\vskip 5mm
\section{An example of type $G$ random variable}

In Maejima and Niiyama \cite{MNi05}, an improper integral 
\begin{equation}\label{4.1}
Z= \int_0^{\infty -} e^{-(B_s+\lambda s)}dS_s
\end{equation}
was studied,
in relation to a stationary solution of the stochastic differential equation
\begin{equation*}
dZ_t = - \lambda Z_{t} dt + Z_{t-} dB_t + dS_t, \quad t\ges 0,
\end{equation*} 
where $\{B_t, t\ges 0\}$ is a standard Brownian motion on $\mathbb R$, $\lambda >0$, and
$\{S _t, t\ges 0\}$ is a symmetric $\alpha$-stable L\'evy process with $0<\alpha \les 2$
on $\mathbb R$, independent of $\{B_t\}$.
They showed that $Z$ is {\it of type $G$} in the sense that $Z$ is a variance mixture
of a standard normal random variable by some infinitely divisible distribution.
Namely, $Z$ is of type $G$ if 
\begin{equation*}
Z\eqd V^{1/2}W
\end{equation*}
for some nonnegative infinitely divisible random variable
$V$ and a standard normal random variable $W$ independent of each other.
Equivalently, $Z$ is of type $G$ if and only if $Z\eqd U_1$, where $\{U_t, t\ges0\}$ 
is given by subordination of a standard Brownian motion.
If $Z$ is of type $G$, then $\mathcal L(V)$ is uniquely 
determined by $\mathcal L(Z)$ (Lemma 3.1 of \cite{S01b}).

The $Z$ in \eqref{4.1} is a special case of those exponential integrals of L\'evy
processes which we are dealing with. Thus Theorem \ref{t3.1} says that 
the law of $Z$ is selfdecomposable. But the class of type $G$ distributions
(the laws of type $G$ random variables) is neither larger nor smaller 
than the class of symmetric selfdecomposable distributions. 
Although the proof that $Z$ is of type $G$ is found in \cite{MNi05},
the research report is not well distributed.
Hence we give their proof below for readers.
We will show that the law of $Z$ belongs to a special 
subclass of selfdecomposable distributions.

\begin{thm}\label{t4.1}
Under the assumptions on $\{B_t\}$ and $\{S_t\}$ stated above, 
$Z$ in \eqref{4.1} is of type $G$ and furthermore the mixing distribution for 
variance, $\mathcal L(V)$, is not only
infinitely divisible but also selfdecomposable.
\end{thm}

\begin{proof}
It is known (Proposition 4.4.4 of 
Dufresne \cite{D90}) that for any $a\in\mathbb R\setminus \{0\}$, $b>0$,
$$
\int_0^{\infty} e^{aB_s-bs} ds \eqd {2}\left (a^2 \Gamma _{2ba^{-2}}\right )^{-1},
$$
where $\Gamma _\gamma$ is the gamma random variable with parameter $\gamma >0$,
namely, $P(\Gamma _\gamma \in B) = \Gamma (\gamma) ^{-1}\int_{B\cap (0,\infty)}
x^{\gamma -1}e^{-x}dx$.
The law of the reciprocal of gamma random variable is
infinitely divisible and, furthermore, selfdecomposable (Halgreen \cite{H79}). 
We have
\begin{align*}
E  \left[ e^{iz Z } \right] 
& = E \left[ \exp \left(iz \int_0^{\infty -} e^{-(B_s+\lambda s)} dS_s 
\right) \right] \\
& = E \left[ E \left[ \left. \exp \left( iz \int_0^{\infty -} 
e^{-(B_s+\lambda s)} dS_s\right)\,\right|\,\{B_s\} \right]\right] ,
\end{align*}
We have $Ee^{izS_t}=\exp(-ct|z|^{\alpha})$ with some $c>0$.  
For any nonrandom measurable function $f(s)$ satisfying $\int_0^{\infty}
|f(s)|^{\alpha}ds<\infty$, we have
$$
E\left [\exp \left( iz \int_0^{\infty -} f(s)dS_s\right)\right ]
= \exp \left( -c|z |^{\alpha}\int_0^{\infty} |f(s)|^{\alpha}ds\right)
$$
(see, e.\,g.\ Samorodnitsky and Taqqu \cite{ST}). Hence
\begin{align*}
E\left [e^{iz Z}\right ]
& = E\left [\exp \left( -c|z |^{\alpha}\int_0^{\infty}
e^{ -\alpha B_s -\alpha \lambda s} ds \right)\right]\\
& =  E\left [\exp\left(-c|z |^{\alpha} 2 
\left (\alpha ^2  \Gamma_{2\alpha ^{-1}\lambda }\right )^{-1}\right)\right ] .
\end{align*}
If we put
$$
H(dx) = P\left (  2 c
\left (\alpha ^2 \Gamma_{2\alpha ^{-1}\lambda  }\right )^{-1}
\in dx\right ),
$$
then 
$$
E[e^{iz Z}] = \int_0^{\infty} e^{-u|z |^{\alpha}} H(du).
$$
This $H$ is the distribution of a positive
infinitely divisible (actually selfdecomposable) random variable.
This shows that $Z$ is a mixture of a symmetric
$\alpha$-stable random variable $S$ with $Ee^{izS}=e^{-|z|^{\alpha}}$ 
in the sense that
\begin{equation}\label{4.2}
Z \eqd \Gamma ^{-1/\alpha}S,
\end{equation}
where $\Gamma$ and $S$ are independent and $\Gamma$ is a gamma random variable with
$\mathcal L(\Gamma^{-1})=H$, that is, 
$\Gamma=(2c)^{-1}\alpha^2 \Gamma_{2\alpha^{-1}\lambda}$.
To see that $Z$ is of type $G$, we need to rewrite \eqref{4.2} as
$$
Z \eqd \Gamma ^{-1/\alpha} S \eqd V^{1/2}W,
$$
for some infinitely divisible random variable $V>0$ independent of a standard 
normal random variable $W$.
Let $S^+_{\alpha/2}$ be a positive strictly $(\alpha/2)$-stable random variable
such that
$$
E\left[\exp(-uS_{\alpha/2}^+)\right]=\exp\left( -(2u)^{\alpha/2}\right),\quad u\ges 0
$$
and $\Gamma$, $W$, and $S_{\alpha/2}^+$ are independent.  Then
$$
S\eqd (S_{\alpha/2}^+)^{1/2} W,
$$
and hence $S$ is of type $G$. Let
$$
V=\Gamma^{-2/\alpha}S_{\alpha/2}^+.
$$
Then
$$
V^{1/2} W=(\Gamma^{-2/\alpha}S_{\alpha/2}^+)^{1/2} W
=\Gamma^{-1/\alpha}(S_{\alpha/2}^+)^{1/2}W
\eqd\Gamma^{-1/\alpha} S\eqd Z.
$$
Using a positive strictly $(\alpha/2)$-stable L\'evy process $\{S_{\alpha/2}^+(t),
t\ges0\}$ independent of $\Gamma$ with $\mathcal L(S_{\alpha/2}^+(1))=S_{\alpha/2}^+$, 
we 
see that 
$$
V\eqd S_{\alpha/2}^+(\Gamma^{-1}).
$$
Since $\Gamma^{-1}$ is selfdecomposable, $V$ is also selfdecomposable due to the
inheritance of selfdecomposability in subordination of strictly stable L\'evy
processes (see \cite{S01b}).
Therefore $Z$ is of type $G$ with $\mathcal L(V)$ being selfdecomposable. Also,
the selfdecomposability of $Z$ again follows.
\end{proof}

In their recent paper \cite{AMR06}, Aoyama, Maejima, and Rosi\'nski
 have introduced a new strict subclass 
(called $M(\mathbb R^d)$) of
the intersection of the class of type $G$ distributions and the class of selfdecomposable
distributions on $\mathbb R^d$ (see Maejima and Rosi\'nski \cite{MR} for the 
definition of 
type $G$ distributions on $\mathbb R^d$ for general $d$).
If we write the polar decomposition of the L\'evy measure $\nu$ by
\begin{equation*}
\nu (B) = \int_K \lambda (d\xi)\int_0^{\infty} 1_B(r\xi)\nu_{\xi}(dr),
\end{equation*}
where $K$ is the unit sphere $\{\xi\in\mathbb R^d\colon |\xi|=1\}$ 
and $\lambda$ is a probability measure on $K$,
then the element of $M(\mathbb R^d)$ is characterized as a symmetric infinitely
divisible distribution such that
\begin{equation*}
\nu_{\xi}(dr) = g_{\xi}(r^2)r^{-1}dr
\end{equation*}
with $g_{\xi}(u)$ being completely monotone as a function of 
$u\in(0,\infty)$ and measurable with respect to $\xi$.
Recall that if we write $\nu_{\xi}(dr) = g_{\xi}(r^2)dr$ instead, this gives a 
characterization of type $G$ distributions on $\mathbb R^d$ (\cite{MR}).
In \cite{AMR06} it is shown that
$$
\{ \text{type $G$ distributions on $\mathbb R$ with selfdecomposable mixing 
distributions}\} \subsetneqq M(\mathbb R).
$$

Now, by Theorem \ref{t4.1} combined with the observation above, we see that
$\mathcal L (Z)$ in \eqref{4.1} belongs to $M(\mathbb R)$.
It is of interest as a
concrete example of random variable whose distribution belongs to $M(\mathbb R)$.

We end the paper with a remark that, by Preposition 3.2 of \cite{CPY01}, 
if $\alpha =2$, our $\mathcal L (Z)$ is also 
Pearson type IV distribution of parameters $\lambda$ and $0$.

\bigskip
{\bf Acknowledgments.}  The authors would like to thank Alexander Lindner and 
Jan Rosi\'nski for their helpful comments while this
paper was written.

\end{document}